\documentclass {article}

\usepackage{latexsym, amsfonts}
\usepackage{amssymb, amsmath}
\usepackage{amsfonts, latexsym}
\usepackage{color}

\newcommand{\im}{\mathrm{im}}
\newcommand{\dom}{\mathrm{dom}}

\newcommand{\rank}{\mathrm{rank}}
\newcommand{\qed}{\hfill $\square$ \vspace{3mm}}

\newtheorem{theorem}{Theorem}
\newtheorem{lemma}[theorem]{Lemma}

\newtheorem{proposition}[theorem]{Proposition}
\newtheorem{remark}[theorem]{Remark}

\date{}

\begin{document}

\title{On semigroups of order-decreasing partial isometries: presentations and ranks of certain ideals}
\author{Gonca Ay\i k$^{1}$,  Hayrullah Ay\i k$^{1}$, V\'itor H. Fernandes$^{2,*}$\vspace{3mm}\\ \small 
$^{1}$Department of Mathematics, \c{C}ukurova University, Adana, Turkey\\ \small 
$^{2}$Center for Mathematics and Applications (NOVA Math) and \\ \small 
Department of Mathematics, NOVA School of Science and Technology,\\ \small 
NOVA University of Lisbon, Lisbon, Portugal}

\date{}
\maketitle
\let\thefootnote\relax\footnote{$^{*}$ Corresponding author.}

\let\thefootnote\relax\footnote{\, \ E-mail adresses:  agonca@cu.edu.tr (Gonca Ay\i k), hayik@cu.edu.tr (Hayrullah Ay\i k), 
vhf@fct.unl.pt (V\'itor H. Fernandes)}	

\let\thefootnote\relax\footnote{\, \ ORCID iD: 0000-0003-2301-0669 (Gonca Ay\i k), 
0000-0001-7569-5597 (Hayrullah Ay\i k), 
0000-0003-1057-4975 (V\'itor H. Fernandes)}

\begin{abstract}
\noindent Let $\mathcal{DP}_{n}^{-}$ denote the monoid of all order-decreasing partial isometries on a chain with $n$ elements, 
and let $\mathcal{ODP}_{n}^{-}$ denote the submonoid of $\mathcal{DP}_{n}^{-}$ consisting of all order-preserving partial isometries. 
In this paper, for each $1\leq r\leq n$, we determine minimal generating sets of the ideals 
$\mathcal{ODP}^{-}(n,r) =\{ \alpha \in \mathcal{ODP}^{-}_{n} : \lvert \im(\alpha) \rvert \leq r \}$ of $\mathcal{ODP}^{-}_{n}$ 
and $\mathcal{DP}^{-}(n,r) =\{ \alpha \in \mathcal{DP}^{-}_{n} : \lvert \im(\alpha) \rvert \leq r \}$ of $\mathcal{DP}^{-}_{n}$, 
thereby determining their ranks. 
Furthermore, we give presentations for the monoids $\mathcal{ODP}^{-}_{n}$ and $\mathcal{DP}^{-}_{n}$.	

\vspace{3mm}
	

\vspace{3mm}


\end{abstract}

\section{Introduction}

For $n\in \mathbb{N}$, let $\mathcal{PT}_{n}$ and $\mathcal{I}_{n}$ be the partial and partial injective transformation semigroups on the finite chain $X_{n}=\{1, \ldots ,n\}$ under its natural order, respectively. A transformation $\alpha\in \mathcal{PT}_{n}$ is called \emph{isometry} or \emph{distance-preserving} if $\lvert x\alpha -y\alpha \rvert=\lvert x-y\rvert$ for all $x,y \in \dom (\alpha)$. If we denote the subset of $\mathcal{PT}_{n}$ consisting of all distance-preserving partial transformations by $\mathcal{DP}_{n}$, then it is clear that $\mathcal{DP}_{n}$ is a submonoid of $\mathcal{I}_{n}$. Since distance-preserving transformations are the primary subject of this study, we proceed with definitions concerning partial injective transformations. A transformation $\alpha\in \mathcal{I}_{n}$ is called \emph{order-preserving} (\emph{order-reversing}) if $x\leq y$ implies $x\alpha\leq y\alpha$ ($y\alpha\leq x\alpha$) and \emph{order-decreasing} if $x\alpha \leq x$ for all $x\in \dom (\alpha)$. A transformation is called monotone if it is order-preserving or order-reversing. Let $\mathcal{I}_{n}^{-}$,\, $\mathcal{IO}_{n}$ and $\mathcal{IM}_{n}$ denote the subsemigroups of $\mathcal{I}_{n}$ consisting of all order-decreasing, order-preserving and monotone partial injective transformations, respectively. It is straightforward to verify that every $\alpha\in \mathcal{DP}_{n}$ is monotone, i.e. $\mathcal{DP}_{n}$ is a subsemigroup of $\mathcal{IM}_{n}$. We also define the following transformations semigroups, which we will concern, 
$$\begin{array}{lcl}
	\mathcal{IM}_{n}^{-} =\mathcal{IM}_{n} \cap \mathcal{I}_{n}^{-}, & & \mathcal{DP}_{n}^{-} =\mathcal{DP}_{n} \cap \mathcal{I}_{n}^{-}, \vspace{1mm}\\
	\mathcal{ODP}_{n} =\mathcal{DP}_{n} \cap \mathcal{IO}_{n} &\mbox{and}& \mathcal{ODP}_{n}^{-} =\mathcal{ODP}_{n} \cap \mathcal{I}_{n}^{-}.
\end{array}$$
Moreover, for each $0\leq r\leq n$, let
$$\mathcal{M}(n,r)=\{ \alpha\in \mathcal{M}_{n} : \lvert \im(\alpha) \rvert \leq r\},$$
where $\mathcal{M}_{n}\in \{ \mathcal{DP}_{n}, \mathcal{ODP}_{n}, \mathcal{IM}_{n}^{-}, \mathcal{DP}_{n}^{-}, \mathcal{ODP}_{n}^{-} \}$. It is clear that 
$\mathcal{M}(n,r)$ is an ideal of $\mathcal{M}_{n}$ for $0\leq r\leq n$. 

\smallskip 

Let $Y$ be any subset of $X_{n}$ and let $1_{Y}$ denote the partial identity of $Y$, i.e. the mapping $1_{Y} : Y\rightarrow Y$ defined by $x1_{Y}=x$ for all $x\in Y$. In particular, 
we use $0_{n}$ and $1_{n}$ to denote the zero (empty) transformation $1_{\emptyset}$ on $X_{n}$ and the identity transformation $1_{X_{n}}$ on $X_{n}$, respectively.  
The \emph{rank} of $\alpha\in \mathcal{I}_{n}$ is defined by $\rank(\alpha) =\lvert \im(\alpha) \rvert$. 

\smallskip 

Let $S$ be an arbitrary semigroup. We denote the subsemigroup generated by a nonempty subset $A$ of $S$, i.e. the smallest subsemigroup of $S$ containing $A$, by $\langle A\, \rangle$. Moreover, if $S=\langle A\, \rangle$, then $A$ is called a \emph{generating set} of $S$. The \emph{rank} of a (finite) semigroup $S$ is defined by $\rank(S) =\min\{ \lvert A\rvert : S=\langle A\, \rangle \}$, and a generating set of $S$ with size $\rank(S)$ is called a \emph{minimal generating set}. An element $a$ of a semigroup $S$ is called an \emph{idempotent} if $a^{2}=a$, and the set of all idempotents of a subset $T$ of $S$ is denoted by $E(T)$. Then notice that $\alpha \in \mathcal{I}_{n}$ is an idempotent if and only if $\alpha =1_{Y}$ for some $Y\subseteq X_{n}$, and thus we get $E(\mathcal{I}_{n}) =\{ 1_{Y} : Y\subseteq X_{n}\}$. An element $a$ of a semigroup $S$ is called \emph{undecomposable} in $S$ if $a\neq st$ for all $s,t\in S\setminus \{a\}$, i.e. $a\in S$ is undecomposable in $S$ if and only if $S\setminus \{a\}$ is a subsemigroup of $S$. It is clear that every generating set of a semigroup contains all undecomposable elements and if the set of all undecomposable elements is a generating set, then it is the only minimal generating set. 
For other terms in semigroup theory, we refer the readers to \cite{GM,H}. 

\smallskip 

We now recall some notions related to the concept of a monoid presentation. For a set $A$, let $A^{*}$ denote the free monoid on $A$ consisting of all finite words over $A$. The empty word is denoted by $\emptyset$. A \emph{monoid presentation} is an ordered pair $\langle A\mid R\rangle$ where $A$ is an alphabet and $R$ is a subset of $A^{*}\times A^{*}$. Each element $(u,v)$ of $R$ is called a \textit{$($defining$)$ relation}, and it is usually written by $u=v$. A monoid $M$ is said to be \textit{defined by a presentation} $\langle A\mid R\rangle$ if $M$ is isomorphic to $A^{*}/\sim_{R}$, where $\sim_R$ denotes the congruence on $A^*$ generated by $R$, i.e. $\sim_{R}$ is the smallest congruence on $A^{*}$ containing $R$. Let  $X$ be a generating set of a monoid $M$ and let $\phi : A\rightarrow M$ be an injective mapping such that $A\phi =X$. If $\Phi : A^{*}\rightarrow M$ is the (unique) homomorphism that extends $\phi$ to $A^{*}$, then we say a relation $u=v$ of $A^{*}$ holds in $M$ via $\Phi$ if $u\Phi =v\Phi$. For more details, see \cite{Lallement:1979} or \cite{Ruskuc:1995}.

\smallskip 

A well-known direct method to obtain a presentation for a monoid is given by the following result; see, for example, \cite[Proposition 1.2.3]{Ruskuc:1995}.

\begin{proposition}\label{provingpresentation}
Let $M$ be a monoid generated by a set $X$ and let $A$ be an alphabet with the same size as $X$ and $R\subseteq A^{*}\times A^{*}$.   
Let $\Phi:A^{*}\rightarrow M$ the homomorphism such that $A\Phi=X$. 
Then $\langle A\mid R\rangle$ is a presentation for $M$ if and only if the following two conditions are satisfied:
\begin{enumerate}
\item all relations from $R$ hold in $M$ via $\Phi$; 
\item if $w,w'\in A^*$ are any two words such that the relation $w=w'$ holds in $M$ via $\Phi$, then $w\sim_R w'$.
\end{enumerate}
\end{proposition}

For a finite monoid, a usual method to find a presentation is described by the following result (adapted to the monoid case from \cite[Proposition 3.2.2] {Ruskuc:1995}). 

\begin{proposition}[Guess and Prove method] \label{ruskuc} 
Let $M$ be a finite monoid generated by a set $X$, let $A$ be an alphabet with the same size as $X$ 
and let $\Phi:A^{*}\rightarrow M$ the homomorphism such that $A\Phi=X$. 
Let $R\subseteq A^{*}\times A^{*}$ and $W\subseteq A^{*}$. Suppose that the following conditions are satisfied:
\begin{enumerate}
\item all relations from $R$ hold in $M$ via $\Phi$; 
\item for each word $w\in A^{*}$, there exists a word $w'\in W$ such that $w\sim_R w'$;
\item $|W|\leqslant |M|$.
\end{enumerate}
Then, $M$ is defined by the presentation $\langle A\mid R\rangle$. 
\end{proposition}

Notice that, if $W$ satisfies the above conditions then, in fact, $\lvert W\rvert =\lvert M\rvert$. Moreover, we say that $W\subseteq A^{*}$ is a set of  \emph{canonical forms} for a finite monoid $M$ with respect to $R\subseteq A^{*} \times A^{*}$, if the conditions of Proposition \ref{ruskuc} are satisfied. 

\smallskip 

The cardinalities of $\mathcal{ODP}_{n}$ and $\mathcal{DP}_{n}$ are determined in \cite{AKU1}, and the cardinalities of $\mathcal{ODP}_{n}^{-}$ and $\mathcal{DP}_{n}^{-}$ are determined in \cite{AKU3}. The ranks of  $\mathcal{ODP}_{n}$ and $\mathcal{DP}_{n}$ are determined in \cite{AKU2}, and the ranks of $\mathcal{ODP}_{n}^{-}$ and $\mathcal{DP}_{n}^{-}$ are determined in \cite{KMU}. Moreover, for each $2\leq r\leq n-1$, the ranks of the ideals $\mathcal{ODP}(n,r)$ and $\mathcal{DP}(n,r)$ of $\mathcal{ODP}_{n}$ and $\mathcal{DP}_{n}$, respectively, are determined in \cite{BYA}. For $1\leq r\leq n$, the subsemigroup $\{ \alpha\in \mathcal{ODP}_{n} : \im(\alpha) \subseteq X_{r} \}$ is studied in \cite{SS}. Presentations for the monoids 
$\mathcal{ODP}_{n}$ and $\mathcal{DP}_{n}$ are given in \cite{FQ1}. Moreover, the monoids of partial isometries of a finite star graph and of a finite cycle graph are studied in \cite{FP,FQ2}. The goal of this study is determine the ranks of the ideals $\mathcal{ODP}^{-}(n,r)$ of $\mathcal{ODP}^{-}_{n}$ and $\mathcal{DP}^{-}(n,r)$ of $\mathcal{DP}^{-}_{n}$. 
Furthermore, we give presentations for the monoids $\mathcal{ODP}^{-}_{n}$ and $\mathcal{DP}^{-}_{n}$.

\section{Generating sets of $\mathcal{ODP}^{-}(n,r)$}

We start with some obvious observations in the following remark which we frequently use explicitly or implicitly.

\begin{remark}\label{r6}
Let $\alpha\in \mathcal{ODP}^{-}_{n}\setminus\{0_n\}$. Then there exists $0\leq k <\min(\dom(\alpha))$ such that $x\alpha =x-k$  for all $x\in \dom(\alpha)$. Furthermore, $\alpha\in E(\mathcal{ODP}^{-}_{n})$ if and only if $k=0$. 
\end{remark}

For each $0\leq r\leq n$,  let  
$$
	J_{r}=\{ \alpha\in \mathcal{ODP}^{-}_{n} : \lvert \im(\alpha) \rvert=r\}
$$ 
and 
$$
 E_{r}=E(J_{r})=\{ 1_{Y} : Y\subseteq X_{n} \mbox{ and }  \lvert Y\rvert =r \}. 
$$
Next let $\pi =\left( \begin{matrix}
	2&3&\cdots&n\\
	1&2&\cdots&n-1
\end{matrix} \right)$. Notice that $\pi\in J_{n-1}$,  $\pi^n=0_n$ and 
$\pi^{k} =\left( \begin{matrix}
	k+1&k+2&\cdots& n\\
	1 & 2 &\cdots&n-k
\end{matrix} \right)$ 
for $1\leq k\leq n-1$.

\begin{lemma}\label{p7}
$J_{r} \subseteq \langle J_{r+1}\rangle$ for $0\leq r\leq n-2$. 
\end{lemma}

\begin{proof}  
The result is clear for $r=0$. So suppose that $r\geq1$. 

\smallskip
	
Let $\alpha\in J_{r}$ and let $0\leq k <\min(\dom(\alpha))$ be such that $x\alpha =x-k$ for all $x\in \dom(\alpha)$. First, let us consider that $\dom(\alpha)=[n-r+1,n]$. Then, noticing that $n-r+1\geq3$, we have $1_{\{1\}\cup \dom(\alpha)},\, 1_{\{2\}\cup \dom(\alpha)} \pi^{k}\in J_{r+1}$ and 
$$\alpha = 1_{\{1\}\cup\dom(\alpha)}\, (1_{\{2\}\cup\dom(\alpha)} \pi^k) \in \langle J_{r+1}\rangle. $$
Next suppose that $\dom(\alpha)\neq [n-r+1,n]$. Then there exists $a\in [n-r+1,n]$ such that $a\not\in\dom(\alpha)$. Let $\beta$ be the partial mapping defined by $\dom(\beta)= \dom(\alpha) \cup\{a\}$ and $x\beta=x-k$ for all $x\in \dom(\beta)$. Then clearly $\beta\in J_{r+1}$. Take $b\in X_{n} \setminus \im(\beta)$. Thus we have
$$\alpha =\beta \hspace{1pt} 1_{\{b\}\cup\im(\alpha)}  \in \langle J_{r+1}\rangle,$$
as required. \qed
\end{proof} 

The next lemma is straightforward to verify.

\begin{lemma}\label{p9}
Let $\gamma\in J_{r}$ and let $\xi\in E_{r}$. Then we have  
\begin{itemize}
	\item [$(i)$] $\gamma \xi\in  J_{r}$ if and only if $\gamma \xi= \gamma$, and
	\item [$(ii)$] $\xi \gamma\in  J_{r}$ if and only if $\xi \gamma= \gamma$. 
\end{itemize}
\end{lemma}

Let $1\leq r\leq n-1$ and let $Y$ be a subset of $[2,n]$ with size $r$. Define the partial mapping $\xi_{Y}$ with domain $Y$ by $x\xi_{Y} =x-1$ for all $x\in Y$. Then, clearly, $\xi_{Y}\in J_{r}$. Moreover, for each $1\leq r\leq n-1$, put
$$F_{r}=\{ \xi_{Y} : Y\subseteq [2,n] \mbox{ and }\lvert Y\rvert=r\}. $$
Notice that we have $\lvert F_{r}\rvert =\binom{n-1}{r}$, $F_{n-1}=\{ \pi\}$ and $J_{n-1}= E_{n-1}\cup F_{n-1}$.    

\begin{lemma}\label{l10}
Let $1\leq r\leq n-1$. Then each element of $E_{r}\cup F_{r}$ is undecomposable in $\mathcal{ODP}^{-}(n,r)$. In particular, each element of $J_{n-1} \cup \{1_{n}\}$ is undecomposable in $\mathcal{ODP}^{-}_{n}$. 
\end{lemma}

\begin{proof}  
Let $\alpha\in E_{r}\cup F_{r}$.  Suppose that $\alpha =\beta \gamma$ for some $\beta, \gamma\in \mathcal{ODP}^{-}(n,r)$. Since $\dom(\alpha) \subseteq \dom(\beta)$ and $\im(\alpha) \subseteq \im(\gamma)$, it is clear that $\beta, \gamma\in J_{r}$. Take $p=\min(\dom(\alpha))$. Suppose that $\beta,\gamma\not\in E_r$. Then, by Remark \ref{r6}, $p\alpha = (p\beta) \gamma < p\beta < p$, whence $p-2\geq p\alpha \in\{p-1,p\}$, which is a contradiction. Thus either $\beta\in E_{r}$ or $\gamma\in E_{r}$. Hence, by Lemma \ref{p9}, we have either $\alpha =\beta$ or $\alpha =\gamma$, and so $\alpha$ is undecomposable in $\mathcal{ODP}^{-}(n,r)$.

\smallskip

Since $J_{n}= \{1_{n}\}$, it is clear that each element of $J_{n-1} \cup \{1_{n}\}$ is undecomposable also in $\mathcal{ODP}^{-}_{n}$.\qed 
\end{proof} 

Now we state and prove the main theorem of this section.

\begin{theorem}\label{t11}
Let $1\leq r\leq n-1$. Then $E_{r}\cup F_{r}$ is the minimal generating set of $\mathcal{ODP}^{-}(n,r)$, and so
$$\rank(\mathcal{ODP}^{-}(n,r)) =\binom{n}{r}+ \binom{n-1}{r}= \frac{2n-r}{n} \binom{n}{r}.$$ 
In particular, $J_{n-1} \cup \{1_{n}\}$ is the minimal generating set of $\mathcal{ODP}^{-}_{n}$, and so $\rank(\mathcal{ODP}^{-}_{n}) =n+2$. 
\end{theorem}

\begin{proof} 
Let $1\leq r\leq n-1$. By Lemma \ref{p7}, it is enough to show that $J_{r} \subseteq \langle E_{r}\cup F_{r} \rangle$. Let $\alpha\in  J_{r} \setminus (E_{r}\cup F_{r})$ with $\dom(\alpha) =Y$ and let $p=\min(Y)$. Then it is clear that $p\geq 3$ and $p\alpha =p-k$ for some $2\leq k\leq p-1$. For each $1\leq i\leq k-1$, if we let $Y_{i}=\{ x-i : x\in Y\}$, then we obtain $\alpha =\xi_{Y} \xi_{Y_{1}} \cdots \xi_{Y_{k-1}} \in \langle F_{r}\rangle$, as required. \qed
\end{proof}

\section{Generating sets of $\mathcal{DP}^{-}(n,r)$}

Notice that $\alpha \in \mathcal{DP}^{-}_{n}$ is both order-preserving and order-reversing if and only if $\lvert \im(\alpha) \rvert \leq 1$. Thus, for all $\alpha\in \mathcal{DP}^{-}_{n} \setminus \mathcal{ODP}^{-}_{n}$, we have $\lvert \im(\alpha) \rvert\geq 2$.
	
\smallskip

For a real number $x$, let $\left\lfloor x\right\rfloor$ denote the greatest integer less than or equal to $x$. 
For $p,q\in X_{n}$,  let 
$$
[p,q]=\{ x\in X_{n} : p\leq x\leq q\},
$$ 
an interval of $X_{n}$, and so $[p,q]=\emptyset$ if $q<p$. 

\smallskip 

Next we examine the elements of $\mathcal{DP}^{-}_{n} \setminus \mathcal{ODP}^{-}_{n}$ a bit more closely.

\begin{lemma}\label{l3}
Let $\alpha\in \mathcal{DP}^{-}_{n} \setminus \mathcal{ODP}^{-}_{n}$ with $\min(\dom(\alpha))=p$. Then we have
\begin{itemize}
	\item [$(i)$] $2\leq p\leq n-1$,\, $\dom(\alpha) \subseteq [p,n]$,\, $\im(\alpha) \subseteq [1,p]$, and
	\item [$(ii)$] $2\leq \lvert \im(\alpha) \rvert \leq \left\lfloor \frac{n+1}{2} \right\rfloor$.
\end{itemize}
\end{lemma}

\begin{proof} 
Since $(i)$ is clear, we prove $(ii)$. Then $2\leq \lvert \im(\alpha) \rvert\leq p\alpha \leq p\leq \left\lfloor \frac{n+1}{2} \right\rfloor$, if $2\leq p\leq \left\lfloor \frac{n+1}{2} \right\rfloor$, and $2\leq \lvert \im(\alpha) \rvert =\lvert \dom(\alpha) \rvert\leq n-p+1 \leq \left\lfloor \frac{n+1}{2} \right\rfloor$, if $\left\lfloor \frac{n+1}{2} \right\rfloor <p\leq n-2$. Thus, in both cases, we obtain $2\leq \lvert \im(\alpha) \rvert \leq \left\lfloor \frac{n+1}{2} \right\rfloor$. \qed
\end{proof}

Since $\mathcal{DP}^{-}(n,1) =\mathcal{ODP}^{-}(n,1)$, we suppose that $2\leq r\leq n$. Then we begin to define some elements of $\mathcal{DP}^{-}(n,r)\setminus\mathcal{ODP}^{-}_{n}$ that, together with $\mathcal{ODP}^{-}(n,r)$, we will prove to be essential for generating $\mathcal{DP}^{-}(n,r)$.

\smallskip

First, for each $2\leq i\leq n-1$, define
\begin{equation*}
\xi_{i}=\left\{ \begin{array}{ll}
	\left( \begin{matrix}
		i&i+1&\cdots &2i-1\\
		i&i-1&\cdots &  1  \end{matrix}\right) &\mbox{if } 2\leq i\leq \left\lfloor \frac{n+1}{2} \right\rfloor  \vspace{2mm}\\
	\left( \begin{array}{cccc}
		i&i+1&\cdots &  n  \\
		i&i-1&\cdots &2i-n \end{array}\right) & \mbox{if }\left\lfloor\frac{n+1}{2}\right\rfloor <i \leq n-1. 
	\end{array}\right.
\end{equation*}
Our next lemma shows the relevance of the transformations we have just defined. 

\begin{lemma}\label{dpger}
Let $\alpha\in\mathcal{DP}^{-}_n\setminus\mathcal{ODP}^{-}_{n}$. Then $\alpha=\xi_i\beta$ with $i=\min(\dom(\alpha))$ and for some $\beta\in \mathcal{ODP}^{-}_{n}$ such that $\rank(\beta)=\rank(\alpha)$. 
\end{lemma}

\begin{proof}
First, notice that $\dom(\alpha) \subseteq \dom(\xi_{i})$ and that $i=\min(\dom(\alpha))$ implies $2\leq i\leq n-1$. Let $0\leq k\leq i-2$ be such that $i\alpha=i-k$. Then $x\alpha=2i-x-k$ for all $x\in\dom(\alpha)$. Define a partial mapping $\beta$ with $\dom(\beta)=\{2i-x : x\in\dom(\alpha)\}$ by $y\beta=y-k$ for all $y\in\dom(\beta)$. Therefore, it is a routine matter to show that $\beta\in \mathcal{ODP}^{-}_{n}$,\, $\rank(\beta)= \rank(\alpha)$ and $\alpha=\xi_i\beta$, as required. \qed
\end{proof}

Observe that 
\begin{equation*}
\rank(\xi_{i})=\left\{ \begin{array}{ll}
	i & \mbox{if } 2\leq i\leq \left\lfloor\frac{n+1}{2}\right\rfloor  \vspace{1mm}\\ 
	n-i+1 & \mbox{if }\left\lfloor\frac{n+1}{2}\right\rfloor <i \leq n-1.
\end{array}\right.
\end{equation*}
Hence, for $2\leq r< \left\lfloor \frac{n+1}{2} \right\rfloor$, we have
\begin{itemize}
	\item[$(i)$] $\rank(\xi_{i})\leq r$ if $2\leq i\leq r$ or $n-r+1\leq i\leq n-1$, and  
	\item[$(ii)$] $\rank(\xi_{i}) >r$ if $r+1\leq i\leq n-r$. 
\end{itemize}
By contrast, we have $\rank(\xi_{i})\leq r$ for all $\left\lfloor \frac{n+1}{2} \right\rfloor \leq r \leq n-1$. 

\smallskip 

Let $2\leq r\leq n-1$.  Define  
\begin{equation*}
G_{r}=\left\{ \begin{array}{ll}
	\{\xi_2,\ldots,\xi_{r},\xi_{n-r+1},\ldots,\xi_{n-1}\} & \mbox{if } 2\leq r < \left\lfloor \frac{n+1}{2} \right\rfloor \vspace{1mm} \\
	\{\xi_2,\ldots,\xi_{n-1}\}& \mbox{if }\left\lfloor \frac{n+1}{2} \right\rfloor \leq r \leq n-1.
\end{array}\right.
\end{equation*}
Thus $G_{r}\subseteq \mathcal{DP}^{-}(n,r)$ and 
$$\lvert G_{r} \rvert= \left\{ \begin{array}{ll}
	2r-2 & \mbox{if }2\leq r < \left\lfloor\frac{n+1}{2}\right\rfloor \vspace{1mm}\\
	 n-2 & \mbox{if }\left\lfloor\frac{n+1}{2}\right\rfloor \leq r \leq n-1.
\end{array}\right.$$
Moreover, it was shown in \cite[Proposition 3.2]{AAIK} that each element of $G_r$ is undecomposable in $\mathcal{IM}^{-}(n,r)$. Since $\mathcal{DP}^{-}(n,r)$ is a subsemigroup of $\mathcal{IM}^{-}(n,r)$, we immediately obtain the following result. 

\begin{lemma}\label{lp12}
Let $2\leq r\leq n-1$. Then each element of $G_{r}$ is undecomposable in $\mathcal{DP}^{-}(n,r)$. 
\end{lemma}

Let $2\leq r\leq n-1$ and define 
$$H_{r} =\{1_{R} \xi_{i} : r+1\leq i\leq n-r,\,\,\, R\subseteq \dom(\xi_{i}),\,\,\, i\in R\,\mbox{ and }\, \lvert R\rvert =r\}. $$
Notice that $H_{r}=\emptyset$ whenever $\left\lfloor\frac{n+1}{2}\right\rfloor \leq r\leq n-1$. On the other hand, since $\rank(\xi_{i})\geq r+1$ for $r+1\leq i\leq n-r$, we notice that $\emptyset\neq H_r\subseteq \mathcal{DP}^{-}(n,r)$ whenever $2\leq r < \left\lfloor\frac{n+1}{2}\right\rfloor$. 

\smallskip 
 
As a consequence of Lemma \ref{dpger}, we have the following proposition:

\begin{proposition}\label{dprrger}
For $2\leq r\leq n-1$, the semigroup $\mathcal{DP}^{-}(n,r)$ is generated by $\mathcal{ODP}^{-}(n,r)\cup G_{r}\cup H_{r}$. 
\end{proposition}

\begin{proof}
Let $\alpha\in\mathcal{DP}^{-}(n,r)\setminus\mathcal{ODP}^{-}_{n}$ and let $i=\min(\dom(\alpha))$. Then, by Lemma \ref{dpger}, $\alpha=\xi_{i}\beta$ for some $\beta\in \mathcal{ODP}^{-}(n,r)$. If either $\left\lfloor \frac{n+1}{2} \right\rfloor \leq r\leq n-1$ or $2\leq r< \left\lfloor \frac{n+1}{2} \right\rfloor$ with $2\leq i\leq r$ or $n-r+1 \leq i\leq n-1$, then $\xi_{i}\in G_{r}$, and so $\alpha$ belongs to the semigroup generated by $\mathcal{ODP}^{-}(n,r) \cup G_{r}$. Thus suppose that $2\leq r< \left\lfloor \frac{n+1}{2} \right\rfloor$ and $r+1\leq i\leq n-r$. Then $\dom(\alpha) \subseteq \dom(\xi_{i})$ and $\lvert \dom(\xi_{i}) \rvert\geq r+1$, and so there exists $\dom(\alpha)\subseteq R\subseteq \dom(\xi_{i})$ such that $\lvert R\rvert=r$. Therefore, $1_{R}\xi_{i} \in H_{r}$ and $\alpha=(1_{R}\xi_{i})\beta$, whence $\alpha$ belongs to the semigroup generated by $\mathcal{ODP}^{-}(n,r) \cup H_{r}$, which completes the proof.  \qed
\end{proof} 

\begin{lemma}\label{lp12+}
Let $2\leq r\leq n-1$. Then each element of $E_{r}\cup F_{r}\cup H_{r}$ is undecomposable in $\mathcal{DP}^{-}(n,r)$. 
\end{lemma}

\begin{proof}
Let $\alpha ,\beta, \gamma\in \mathcal{DP}^{-}(n,r)$ be such that $\alpha =\beta \gamma$ and $\rank(\alpha)=r$. Then $\dom(\alpha) =\dom(\beta)$, $\im(\beta)= \dom(\gamma)$ and $\im(\gamma) =\im(\alpha)$. 

\smallskip 

First, let us suppose that $\alpha\in E_{r}\cup F_{r}$. Then either $\beta,\gamma\in \mathcal{ODP}^{-}(n,r)$ or $\beta, \gamma\in \mathcal{DP}^{-}(n,r) \setminus \mathcal{ODP}^{-}_{n}$. By contradiction, suppose that $\beta, \gamma\in \mathcal{DP}^{-}(n,r) \setminus \mathcal{ODP}^{-}_{n}$. Let $p=\min(\dom(\alpha))$ and $q=\max(\dom(\alpha))$. Then $q= \max(\dom(\beta))$ and $q\beta= \min(\im(\beta))= \min(\dom(\gamma)$. Let $k,\ell \geq0$ be such that $p\beta =p-k$ and $(q\beta) \gamma= q\beta-\ell$. Then $q\beta=2p-k-q$, and so $q\alpha =q(\beta \gamma) = q\beta -\ell= 2p-k-q-\ell$, whence $2p-q-k-\ell\in\{q,q-1\}$, which is a contradiction since $2p-2q\leq -2$ and $k+\ell \geq 0$. Thus $\beta, \gamma\in \mathcal{ODP}^{-}(n,r)$, and so, by Lemma \ref{l10}, we have $\beta=\alpha$ or $\gamma=\alpha$. Therefore, $\alpha$ is undecomposable in $\mathcal{DP}^{-}(n,r)$. 

\smallskip 

Now suppose $\alpha\in H_{r}$. Let $i=\min(\dom(\alpha))$. Then $i=i\alpha =i(\beta \gamma) \leq i\beta \leq i$, whence $i\beta=i$. If $\beta\in \mathcal{ODP}^{-}_n$, then $\beta=1_{\dom(\beta)}$, and so $\gamma =\beta \gamma= \alpha$. On the other hand, if $\beta \not\in \mathcal{ODP}^{-}_{n}$, since $i\alpha =i= i\beta$ and $\dom(\alpha) =\dom(\beta)$, we get $\beta=\alpha$. Therefore, $\alpha$ is undecomposable in $\mathcal{DP}^{-}(n,r)$. \qed 
\end{proof}

Recall the well-known binomial identity 
$\sum\limits_{k=p}^{\ell} \binom{k}{p} =\binom{\ell+1}{p+1} \quad \mbox{for }\ell\geq p\geq 0$,
(called hockey-stick identity, Christmas tree identity or summation on the upper index; for example, see \cite{GKP}), 
from which it is easy to deduce the following two binomial identities: 
\begin{equation}\label{hs}
\begin{array}{ll}
	\sum\limits_{k=p+1}^{\ell}\binom{k-1}{p-1}=\binom{\ell}{p}-1 & \mbox{for }\ell\geq p\geq 1,\vspace{1mm}\\ 
	\sum\limits_{k=m+1}^{\ell-p}\binom{\ell-k}{p-1}=\binom{\ell-m}{p}-1 & \mbox{for }1\leq p\leq \ell-m. 
\end{array}
\end{equation}

\begin{lemma}\label{hr}
Let $m=\left\lfloor \frac{n+1}{2} \right\rfloor$. For $2\leq r<m$,\, $\lvert H_{r}\rvert= \binom{m}{r}+ \binom{n-m}{r} -2$. 
\end{lemma}

\begin{proof} Let $r+1\leq i \leq n-r$. Then 
$$\lvert \dom(\xi_{i}) \setminus \{i\} \rvert= \left\{ \begin{array}{ll}
	i-1 & \mbox{if }r+1\leq i\leq m\\
	n-i & \mbox{if }m < i\leq n-r,
\end{array}\right. $$
whence
$$\lvert \{R : R\subseteq\dom(\xi_i),\,\, i\in R\mbox{ and } \lvert R\rvert=r \}\rvert= \left\{ \begin{array}{ll}
	\binom{i-1}{r-1} & \mbox{if }r+1\leq i\leq m \vspace{1mm}\\
	\binom{n-i}{r-1} & \mbox{if }m <i\leq n-r,
\end{array}\right.$$
and so, by (\ref{hs}),\, 
$\lvert H_{r} \rvert= \sum\limits_{i=r+1}^{m} \binom{i-1}{r-1} + \sum\limits_{i=m+1}^{n-r} \binom{n-i}{r-1} = \binom{m}{r}+\binom{n-m}{r}-2$. \qed 
\end{proof}

From Theorem \ref{t11}, Proposition \ref{dprrger} and Lemmas \ref{lp12}, \ref{lp12+} and \ref{hr}, we can directly derive the main result of this section as follows. 

\begin{theorem}\label{t16}
Let $m=\lfloor \frac{n+1}{2}\rfloor$ and $2\leq r\leq n-1$. Then $E_{r}\cup F_{r}\cup G_r\cup H_r$ is the minimal generating set of $\mathcal{DP}^{-}(n,r)$, and so
$$\rank(\mathcal{DP}^{-}(n,r)) =\!\left\{\!\!\! \begin{array}{ll}
	\frac{2n-r}{n}\binom{n}{r}+ \binom{m}{r}+\binom{n-m}{r} + 2r-4&\!\!\mbox{if }2\leq r < \left\lfloor\frac{n+1}{2}\right\rfloor \vspace{2mm}\\
	\frac{2n-r}{n}\binom{n}{r}+n-2 & \!\!\mbox{if }\left\lfloor\frac{n+1}{2}\right\rfloor \leq r \leq n-1.
\end{array} \right.$$ 
In particular, $J_{n-1}\cup G_{n-1} \cup \{1_{n}\}$ is the minimal generating set of $\mathcal{DP}^{-}_{n}$, and so $\rank(\mathcal{DP}^{-}_{n}) =2n$. 
\end{theorem}

\section{A Presentation for $\mathcal{ODP}^{-}_{n}$} 

In this section, we aim to find a monoid presentation of $\mathcal{ODP}^{-}_{n}$. 

\smallskip

For each $1\leq i\leq n$, let $Y_{i}= X_{n} \setminus \{i\}$. Then $X=\{ 1_{Y_{1}}, \ldots, 1_{Y_{n}}, \pi\}$ is the minimum generating set of the monoid $\mathcal{ODP}^{-}_{n}$. We let $A=\{a_{1}, \ldots, a_{n},b\}$ and define a mapping $\phi :A \rightarrow X$ by $a_{i}\phi =1_{Y_{i}}$ for every $1\leq i\leq n$ and $b\phi =\pi$. Let $\Phi : A^{*} \rightarrow \mathcal{ODP}^{-}_{n}$ be the unique homomorphism such that $\Phi_{\mid_{A}} =\phi$. Then it is a routine matter to check that all the following relations hold in $\mathcal{ODP}^{-}_{n}$ via $\Phi$:

\begin{itemize}
	\item [$R_{1}$]: $a_{i}^{2}=a_{i}$ for $1\leq i\leq n$,
	\item [$R_{2}$]: $a_{i}a_{j}=a_{j}a_{i}$ for $1\leq i<j\leq n$,
	\item [$R_{3}$]: $ba_{n}=b$,  $a_{1}b=b$ and $a_{i}b=ba_{i-1}$ for $2\leq i\leq n$,
	\item [$R_{4}$]: $b^{n}=a_{1} \cdots a_{n}$ and $b^{n+1}=b^{n}$.
\end{itemize}

Denote by $\sim$ the congruence of $A^{*}$ generated by the relations $R_{1}\-- R_{4}$. 

\smallskip

Observe that a simple line of reasoning by induction using the relations $R_{4}$ allows us to deduce that  $b^{k}\, \sim\, a_{1}\cdots a_{n}$ for any $k\geq n$. Moreover, we have the following lemma. 

\begin{lemma}\label{l17} 
Let $1\leq j\leq n-1$. Then
\begin{itemize}
\item[$(i)$] $a_{i} b^{j}\, \sim\, b^{j}$ for $1\leq i \leq j$,
\item[$(ii)$] $a_{j+1}\cdots a_{n}b^{j}\, \sim\, a_{1} \cdots a_{n}$. 
\end{itemize}
\end{lemma} 

\begin{proof}
$(i)$ First, we show by induction on $k$ that $a_{k} b^{k}\, \sim\, b^{k}$ for $1\leq k\leq n-1$. If $k=1$, then $a_{1}b\, \sim_{R_{3}}\, b$. Next let $1\leq k< n-1$ and suppose, by induction hypothesis, that $a_{k} b^{k}\, \sim\, b^{k}$. Then 
$$a_{k+1} b^{k+1} = a_{k+1} b b^{k} \,\sim_{R_{3}}\, b a_{k} b^{k} \,\sim \, b b^{k}=b^{k+1},$$
as required. Hence, we have $a_{i}b^{j}=a_{i}b^{i}b^{j-i} \,\sim\, b^{i}b^{j-i}=b^{j}$ for $1\leq i \leq j$. 
 
\smallskip 

$(ii)$ To start, notice that $a_{1}\cdots a_{n} b \, \sim_{R_{4}}\, b^{n+1} \,\sim_{R_{4}}\, b^{n}  \,\sim_{R_{4}}\, a_{1}\cdots a_{n}$. Now we show by induction on $k$ that $a_{k+1}\cdots a_{n}b^{k} \,\sim\, a_{1} \cdots a_{n}$ for $1\leq k\leq n-1$. If $k=1$, then 
$$
a_{2}\cdots a_{n} b\, \sim_{R_3}\, a_{2}\cdots a_{n} a_{1}b \, \sim_{R_2}\, a_{1}\cdots a_{n} b\, \sim_{R_4}\, a_{1}\cdots a_{n}.
$$ 
Next let $1\leq k< n-1$ and suppose that $a_{k+1}\cdots a_{n}b^{k} \,\sim\, a_{1} \cdots a_{n}$. Then 
\begin{eqnarray*}
a_{k+2}\cdots a_{n}b^{k+1} & \sim & a_{k+2}\cdots a_{n} a_{k+1}b^{k+1} \quad (\mbox{by } (i))\\
& \sim_{R_{2}} & a_{k+1}\cdots a_{n}b^{k}b \\
& \sim &a_{1}\cdots a_{n} b \qquad\mbox{(by the induction hypothesis)}\\ 
& \sim_{R_{4}} & a_{1}\cdots a_{n},
\end{eqnarray*}
which completes the proof. \qed 
\end{proof} 

It was shown in \cite[Theorem 2.4]{AKU3} that $\lvert \mathcal{ODP}^{-}_{n} \rvert =2^{n+1}-n-1$. Furthermore, it is well-known that the monoid presentation $\mathcal{P}_{E}= \langle\, a_{1}, \ldots, a_{n} \mid R_{1}\-- R_{2}\,\rangle$ is a (minimal) monoid presentation of $E(\mathcal{I}_{n})$, which is isomorphic to the free semilattice on $X_{n}$. Clearly 
$$W_{E} = \{ a_{i_{1}}\cdots a_{i_{r}} :1\leq i_{1}< \cdots< i_{r}\leq n\}\cup\{ \emptyset\}$$
is a set of canonical forms for $E(\mathcal{I}_{n})$. Now let $W=W_{E} \cup W_{F} \cup \{b,\ldots, b^{n-1}\}$, where
$$W_{F}=\bigcup_{k=1}^{n-2}\left\{ a_{i_{1}} \cdots a_{i_{r}}b^{k} : k+1\leq i_1<\cdots <i_r\leq n \mbox{ and }1\leq r<n-k \right\}.$$
Then notice that $\lvert W\rvert =2^{n}+ \left(\sum\limits_{k=1}^{n-2} 2^{n-k}-2 \right) +n-1= 2^{n+1}-n-1$.

\begin{theorem}\label{t18}
The monoid $\mathcal{ODP}^{-}_{n}$ is defined by the monoid presentation
$$
\mathcal{P}=\langle\, a_{1}, \ldots, a_{n},\, b \mid R_{1}\--R_{4}\, \rangle .
$$
\end{theorem}

\begin{proof} 
Since all relations $R_{1}\-- R_{4}$ hold in $\mathcal{ODP}^{-}_{n}$ and $\lvert W\rvert= \lvert \mathcal{ODP}^{-}_{n} \rvert$, it remains to show that, for any $w\in A^{*}$, there exists $w'\in W$ such that $w\sim w'$. 

\smallskip

Let $w\in A^{*}$. If $w=\emptyset$, then $w\in W_{E}$. So suppose that $w\neq \emptyset$. If $w=b^{k}$ for some $k\geq 1$, then $w\in W$ for $1\leq k\leq n-1$, and $w\, \sim\, a_{1}\cdots a_{n}\in W_{E}$ for $k\geq n$. If $w$ does not contain the letter $b$, then $w\, \sim\, w'$ for some $w'\in W_{E}$ since $\mathcal{P}_{E}$ is a presentation for $E(\mathcal{I}_{n})$ and $W_{E}$ is a set of canonical forms for $E(\mathcal{I}_{n})$. 

\smallskip

Finally, suppose that $w$ contains the word $b$ and some of the letters $a_{1}, \ldots, a_{n}$. Then, by successively and appropriately applying the relations $R_{3}$, we get $w\, \sim_{R_{3}}\, a_{\ell_{1}}\cdots a_{\ell_{r}} b^{k}$ for some $k\geq 1$,\, $r\geq 1$ and $1\leq \ell_{1},\ldots,\ell_{r}\leq n$. From the previous case, we have $a_{\ell_{1}} \cdots a_{\ell_{r}} \,\sim\, a_{j_{1}}\cdots a_{j_{s}}\in W_{E}$ for some $1\leq s\leq n$ and $1\leq j_{1}<\cdots < j_{s}\leq n$. If $k\geq n$, then 
$$w\, \sim\, a_{j_1}\cdots a_{j_s} b^{k} \,\sim\, a_{j_{1}}\cdots a_{j_{s}} a_{1}\cdots a_{n}\, \sim\, a_{1}\cdots a_{n} \in W_{E},$$
by $R_{1}$ and $R_{2}$. Thus let us suppose that $1\leq k\leq n-1$. Hence, by $R_{2}$ and Lemma \ref{l17} $(i)$, we obtain $w\,\sim\,  a_{i_{1}}\cdots a_{i_{r}} b^{k}$ for some $m+1\leq i_{1} < \cdots < i_{r}\leq n$ and $1\leq r\leq n-k$. If $k=n-1$, then $r=1$,\, $i_{1}=n$ and 
$$w\, \sim\, a_{n} b^{n-1}\,\sim\, a_{1}\cdots a_{n} \in W_{E},$$ 
by Lemma \ref{l17} $(ii)$. If $r=n-k$, then 
$$w\, \sim\, a_{k+1}\cdots a_{n} b^{k}\, \sim\, a_{1}\cdots a_{n}\in W_{E},$$ 
by Lemma \ref{l17} $(ii)$. If $k<n-1$ and $r<n-k$, then  
$$w\, \sim\,  a_{i_{1}}\cdots a_{i_{r}} b^{k} \in W,$$ 
as required. \qed
\end{proof}

\section{A Presentation for $\mathcal{DP}^{-}_{n}$}

Let us consider the minimum generating set 
$$Y=X\cup\{\xi_{2}, \ldots, \xi_{n-1}\} =\{ 1_{Y_{1}}, \ldots, 1_{Y_{n}}, \pi,\xi_{2} ,\ldots, \xi_{n-1}\}$$ 
of the monoid $\mathcal{DP}^{-}_{n}$. 
Let $B=A\cup\{c_{2},\dots,c_{n-2}\}=\{a_{1}, \ldots, a_{n},b,c_{2},\dots,c_{n-2}\}$ and define the mapping $\psi :B \rightarrow Y$ that extends the mapping $\phi : A \rightarrow X$ (i.e. $a_{i}\psi =1_{Y_{i}}$ for  $1\leq i\leq n$ and $b\psi =\pi$) and $c_{i}\psi=\xi_{i}$ for  $2\leq i\leq n-1$. Denote by $\Psi : B^{*} \rightarrow \mathcal{DP}^{-}_{n}$ the unique homomorphism such that $\Psi_{\mid_{B}} =\psi$. Then, obviously, the relations $R_{1}\-- R_{4}$ also hold in $\mathcal{DP}^{-}_{n}$ via $\Psi$ and it is a routine matter to check that all the following relations hold in $\mathcal{DP}^{-}_{n}$ via $\Psi$ too: 
\begin{itemize}
	\item [$R_{5}$]: $c_{i} c_{j} =a_{1}\cdots a_{n}$ for $2\leq i<j\leq n-1$,
	\item [$R_{6}$]: $c_{i} c_{j}=a_{1}\cdots a_{i-1}b^{2(i-j)}a_{j+1}\cdots a_{n}$  for $2\leq j\leq i\leq n-1$,
	\item [$R_{7}$]: $c_{i}^{3} = c_{i}^{2}$ for $2\leq i\leq n-1$,
	\item [$R_{8}$]: $a_{j} c_{i} = c_{i}$ for $1\leq j<i\leq n-1$ and for $3\leq 2i-1 <j\leq n$,
	\item [$R_{9}$]: $c_{i} a_{j} = c_{i}$ for $2\leq i<j\leq n$ and for $1\leq j < 2i-n$ and $\left\lfloor \frac{n+1}{2} \right\rfloor <i\leq n-1$,
	\item [$R_{10}$]: $a_{j} c_{i} = c_{i} a_{2i-j}$ for $2\leq i\leq j\leq \min\{2i-1,n\}$,
	\item [$R_{11}$]: $b c_{n-1} = a_{2}\cdots a_{n-1} b$ and $b c_{i} =c_{i+1} a_{i+2}\cdots a_{n} b$ for $2\leq i\leq n-2$,
	\item [$R_{12}$]: $c_{2} b = a_{3}\cdots a_{n} b$ and  $c_{i} b = a_{2}\cdots a_{i-1} b c_{i-1}$ for $3\leq i\leq n-1$, and 
	\item [$R_{13}$]: $c_{i} b^{2} = c_{i-1} a_{i-1}$ for $4\leq i \leq n-1$. 
\end{itemize}

Note that $\alpha^{0}=1_{n}$ for all $\alpha \in \mathcal{DP}^{-}_{n}$ and that $w^{0}$ is the empty word for all $w\in B^{*}$.

Denote by $\approx$ the congruence of $B^{*}$ generated by the relations $R_{1}\-- R_{13}$. 

\begin{lemma}\label{l1predp} 
For $n\geq 4$, $c_{2} a_{2}\, \approx\, a_{4}\cdots a_{n} b^{2}\, \approx\, c_{3} a_{4} \cdots a_{n} b^{2}$. 
\end{lemma}

\begin{proof}
We start by proving the first relation. 
\begin{eqnarray*}
	c_{2} a_{2}&\sim_{R_{10}}& a_{2} c_{2}\, \sim_{R_{8}}\, a_{2} a_{1} a_{4} \cdots a_{n} c_{2} \, \sim_{R_{2}}\, a_{1} a_{2} a_{4} \cdots a_{n} c_{2}  \, \sim_{R_{2}}\, a_{1}c_{3}^{2} c_{2} \\
	&\sim_{R_{8}}& c_{3}^{2} c_{2}\, \sim_{R_{6}}\, c_{3} a_{2} b^{2} a_{3} \cdots a_{n} \, \approx\, c_{3} b^{2} a_{3} \cdots a_{n} \,\,\, (\mbox{by Lemma \ref{l17} $(i)$})\\ 
	&\sim_{R_{12}}& a_{2} b c_{2} b a_{3} \cdots a_{n} \, \sim_{R_{12}}\, a_{2} b a_{3}\cdots a_{n} b a_{3} \cdots a_{n} \\
	&\sim_{R_{3}}& a_{2} b a_{3}\cdots a_{n}  a_{4} \cdots a_{n} b \, \sim_{R_{1},R_{2}}\, a_{2} b a_{3}\cdots a_{n}  b \,\sim_{R_{3}}\, a_{2} a_{4}\cdots a_{n}  b^{2} \\
	&\sim_{R_{2}}& a_{4}\cdots a_{n}  a_{2} b^{2}\, \sim_{R_{4}}\, a_{4}\cdots a_{n}  b^{2}.
\end{eqnarray*} 
Furthermore, since $a_{2} b a_{3}\cdots a_{n}  b \,\approx\, a_{4}\cdots a_{n}  b^{2}$, we have 
\begin{eqnarray*}
	c_{3} a_{4}\cdots a_{n}  b^{2}&\sim_{R_{9}}& c_{3} b^{2}\, \sim_{R_{12}}\, a_{2} b c_{2} b  \, \sim_{R_{12}}\, a_{2} b a_{3}\cdots a_{n} b\, \approx\, a_{4}\cdots a_{n}  b^{2},
\end{eqnarray*}
as required. \qed
\end{proof}

\begin{lemma}\label{l2predp} 
Let $2\leq i \leq n-1$ and $u\in A^{*}$. Then we have the following: 
\begin{itemize}
	\item[$(i)$] $c_{i} a_{1} \cdots a_{n}\, \approx\, a_{1}\cdots a_{n}$;
	\item[$(ii)$] if $(c_{i} u)\Psi =0_{n}$, then $c_{i} u\, \approx\, a_{1} \cdots a_{n}$; 
	\item[$(iii)$] $c_{i} b^{i-1}\, \approx\, v$ for some $v\in A^{*}$; 
	\item[$(iv)$] if $n\geq 5$ is odd, then $c_{n-1} b^{n-3}\, \approx\, c_{\frac{n+1}{2}}a_3\cdots a_{\frac{n+1}{2}}$; 
	\item[$(v)$] if $1\leq j\leq i$, then $c_{i} a_{1}\cdots a_{j-1} a_{j+1}\cdots a_{n}\, \approx\, v$ for some $v\in A^{*}$; 
	\item[$(vi)$] if the rank of $(c_{i} u)\Psi$ is at most $1$, then $c_{i} u\, \approx\, v$ for some $v\in A^{*}$. 
\end{itemize}  
\end{lemma} 

\begin{proof} $(i)$ We have 
\begin{eqnarray*}
	c_{i} a_{1} \cdots a_{n} &\sim_{R_{2}}& c_{i} a_{1} \cdots a_{i-1} a_{i+1} \cdots a_{n}a_{i}\, \sim_{R_{6}}\, c_{i} c_{i}^{2} a_{i}\, \sim_{R_{7}}\, c_{i}^{2} a_{i}\\ 
	&\sim_{R_{6}}& a_{1} \cdots a_{i-1} a_{i+1} \cdots a_{n} a_{i}\, \sim_{R_5}\, a_{1} \cdots a_{n}. 
\end{eqnarray*} 
	
$(ii)$ First, since $(a_{1} \cdots a_{n})\Phi=0_{n}$, and so $(a_{1} \cdots a_{n} u)\Phi=0_{n}$, by Theorem \ref{t18}, we note that $a_{1} \cdots a_{n} u\, \approx\, a_{1} \cdots a_{n}$ is a consequence of the relations $R_{1}\-- R_{4}$.  
	
\smallskip 
	
Suppose that $(c_{i}u)\Psi =0_{n}$, i.e. $\im(c_{i}\Psi) \cap \dom(u\Psi)= \emptyset$. Since $(a_{1}\cdots a_{i}u)\Phi =u\Phi$ if $i\leq \left\lfloor \frac{n+1}{2} \right\rfloor$, and since $(a_{2i-n}\cdots a_{i} u)\Phi = u\Phi$ if $i > \left\lfloor \frac{n+1}{2} \right\rfloor$, it follows from Theorem \ref{t18} that $a_{1}\cdots a_{i} u\, \approx\, u$ if $i\leq \left\lfloor \frac{n+1}{2} \right\rfloor$, and that $a_{2i-n}\cdots a_{i} u\, \approx\,  u$ if $i> \left\lfloor \frac{n+1}{2} \right\rfloor$. By taking into account $(i)$, we conclude that 
\begin{eqnarray*}
	c_{i} u &\approx& c_{i} a_{1}\cdots a_{i} u \, \sim_{R_{9}}\, c_{i} a_{i+1}\cdots a_{n} a_{1}\cdots a_{i} u \, \sim_{R_{2}}\, c_{i} a_{1}\cdots a_{n} u \\
	&\approx& a_{1}\cdots a_{n}u \,\approx\, a_{1}\cdots a_{n}
\end{eqnarray*} 
if $i\leq \left\lfloor \frac{n+1}{2} \right\rfloor$, and that	
\begin{eqnarray*}
	c_{i} u &\approx& c_{i} a_{2i-n}\cdots a_{i} u \, \sim_{R_{9}}\, c_{i} a_{1} \cdots a_{2i-n-1} a_{i+1}\cdots a_{n} a_{2i-n}\cdots a_{i} u \\ 
	&\sim_{R_{2}}& c_{i} a_{1}\cdots a_{n} u\, \approx\, a_{1}\cdots a_{n} 
\end{eqnarray*}
if $i> \left\lfloor \frac{n+1}{2} \right\rfloor$.
	
\smallskip    
	
$(iii)$ We proceed by induction on $i$. If $i=2$, then $c_{2} b\, \sim_{R_{12}}\, a_{3}\cdots a_{n} b\in A^{*}$. Let $2\leq i< n-2$ and suppose that $c_{i} b^{i-1}\, \approx\, v$ for some $v\in A^{*}$. Then, by the induction hypothesis, we have 
$$ c_{i+1} b^{i} \, \sim_{R_{12}}\, a_{2} \cdots a_{i} b c_{i} b^{i-1}\, \approx\, a_{2} \cdots a_{i} b v \in A^{*},$$
which proves the property. 
	
\smallskip 
	
$(iv)$ If $n=5$, then $c_{4}b^{2}\, \sim_{R_{14}}\, c_{3} a_{3}$. First we prove by induction on $k$ that $c_{n-1} b^{n-3}\, \approx\, c_{n-k} b^{n-2k-1} a_3 \cdots a_{k+1}$ for $2\leq k \leq \frac{n-1}{2}$. Since 
$$c_{n-1}b^{n-3}\, \sim_{R_{13}}\, c_{n-2} a_{n-2} b^{n-5}\, \sim_{R_{3}}\, c_{n-2} b^{n-5} a_3,$$ 
the relation is proved for $k=2$. Let $2\leq k < \frac{n-1}{2}$ and suppose that $c_{n-1} b^{n-3}\, \approx\, c_{n-k} b^{n-2k-1} a_{3} \cdots a_{k+1}$. Then  
\begin{eqnarray*}
	c_{n-1} b^{n-3} &\approx & c_{n-k} b^{n-2k-1} a_{3} \cdots a_{k+1}\\ 
	&\sim_{R_{13}}& c_{n-k-1} a_{n-k-1} b^{n-2k-3} a_{3} \cdots a_{k+1} \\
	& \sim_{R_3} & c_{n-k-1} b^{n-2k-3} a_{k+2} a_3 \cdots a_{k+1} \\ 
	&\sim_{R_2}&  c_{n-k-1} b^{n-2k-3} a_{3} \cdots a_{k+2}, 
\end{eqnarray*} 
which completes the proof by induction. Therefore, in particular for $k=\frac{n-1}{2}$, we have $c_{n-1} b^{n-3}\, \approx\, c_{\frac{n+1}{2}}a_3\cdots a_{\frac{n+1}{2}}$. 
	
\smallskip 
	
$(v)$ If $j=i$, then 
$$c_{i} a_{1}\cdots a_{i-1} a_{i+1}\cdots a_{n}\, \sim_{R_{6}} \, c_{i}c_{i}^{2} \, \sim_{R_{7}}\, c_{i}^{2}\, \sim_{R_{6}}\, a_{1} \cdots a_{i-1}a_{i+1} \cdots a_{n} \in A^{*}.$$ 
If $2\leq j <i$, then 
\begin{eqnarray*} 
	c_{i} a_{1}\cdots a_{j-1} a_{j+1}\cdots a_{n} &\sim_{R_{6}}& c_{i} c_{j}^{2}\, \sim_{R_{6}}\,  a_{2}\cdots a_{i-1} b^{2(i-j)} a_{j+1} \cdots a_{n} c_{j} \\
	&\sim_{R_{8}}& a_{2}\cdots a_{i-1} b^{2(i-j)} a_{j+1} \cdots a_{n} a_{1} \cdots a_{j-1} c_{j} \\
	&\sim_{R_{2}}& a_{2}\cdots a_{i-1} b^{2(i-j)} a_{1}\cdots a_{j-1} a_{j+1}\cdots a_{n}  c_{j} \\
	&\sim_{R_{6}}& a_{2}\cdots a_{i-1} b^{2(i-j)} c_{j}^{3} \, \sim_{R_{7}}\, a_{2}\cdots a_{i-1}b^{2(i-j)} c_{j}^{2} \\
	&\sim_{R_{6}} & a_{2}\cdots a_{i-1}b^{2(i-j)} a_{1} \cdots a_{j-1} a_{j+1} \cdots a_{n} \in A^{*}. 
\end{eqnarray*} 
Finally, let $j=1$ and consider several cases depending on $i$. If $i> \left\lfloor \frac{n+1}{2} \right\rfloor$, then  
$$
c_{i} a_{2} \cdots a_{n}\, \sim_{R_{9}}\, c_{i} a_{1} a_{2}\cdots a_{n}\, \approx\, a_{1} \cdots a_{n} \in A^{*},
$$ 
by $(i)$. 
If $i\leq \left\lfloor \frac{n+1}{2} \right\rfloor$ and $2i-1\leq n-1$, then 
\begin{eqnarray*} 
	c_{i} a_{2} \cdots a_{n} &\sim_{R_{2}}& c_{i} a_{i+1} \cdots a_{n} a_{i}\cdots a_{2}\, \sim_{R_{9}}\, c_{i} a_{i}\cdots a_{2} \, \sim_{R_{10}} \, a_{i} \cdots a_{2i-2} c_{i}\\	
	&\sim_{R_{8},R_{2}}& a_{1}\cdots a_{2i-2} c_{i}\, \sim_{R_{8}}\, a_{1}\cdots a_{2i-2} a_{2i}\cdots a_{n} c_{i}\, \sim_{R_{7}}\, c_{2i-1}^2 c_{i}\\
	&\sim_{R_{6}}& c_{2i-1} a_{2} \cdots a_{2i-2} b^{2(i-1)} a_{i+1} \cdots a_{n}\\  
	&\approx& c_{2i-1} b^{2i-2} a_{i+1} \cdots a_{n} \quad (\mbox{by Lemma \ref{l17} $(i)$})\\ 
	&\approx& v a_{i+1} \cdots a_{n} \in A^{*} \qquad\,\,\, (\mbox{by } (iii)) 
\end{eqnarray*} 
for some $v\in A^{*}$. At last, suppose that $i\leq \left\lfloor \frac{n+1}{2} \right\rfloor$ and $2i-1=n$. So $n$ must be odd and $i=\frac{n+1}{2}$. If $n=3$, then $i=2$ and, by Lemma \ref{l1predp}, $c_{2} a_{2} a_{3}\, \approx\, a_{4} \cdots a_{n} b^{2} a_{3} \in A^{*}$. Suppose that $n\geq 5$. Then
\begin{eqnarray*}
	c_{\frac{n+1}{2}}a_{2}\cdots a_{n} \!& \sim_{R_{1},R_{2}} & 
	c_{\frac{n+1}{2}}a_{3}\cdots a_{\frac{n+1}{2}} a_{2} a_{3} a_{\frac{n+1}{2}} \cdots a_{n}\\
	&\approx& c_{n-1}b^{n-3} a_{2} a_{3} a_{\frac{n+1}{2}} \cdots a_{n} \qquad {(\mbox{by $(iv)$})}\\
	&\approx&c_{n-1}a_{1} \cdots a_{n-3}b^{n-3} a_{2} a_{3}  a_{\frac{n+1}{2}} \cdots a_{n}\,\,\,  {(\mbox{by Lemma \ref{l17} $(i)$})}\\
	&\sim_{R_{3}}&c_{n-1}a_{1} \cdots a_{n-3} a_{n-1} a_{n} b^{n-3}  a_{\frac{n+1}{2}} \cdots a_{n}\\
	&\sim_{R_{1},R_{2}}&c_{n-1}(a_{1} \cdots a_{n-3} a_{n-1} a_{n})(a_{1} \cdots a_{n-3} b^{n-3}  a_{\frac{n+1}{2}} \cdots a_{n})\\
	&\sim_{R_{6}}& c_{n-1} c_{n-2}^{2} c_{n-2} c_{\frac{n-1}{2}} \, \sim_{R_{7}}\, (c_{n-1} c_{n-2}) (c_{n-2} c_{\frac{n-1}{2}})\\
	&\sim_{R_{6}}&\!(a_{1}\cdots a_{n-2}b^{2}a_{n-1}a_{n})(a_{1} \cdots a_{n-3} b^{n-3} a_{\frac{n+1}{2}} \cdots a_{n})\in A^{*}.
\end{eqnarray*}
	
\smallskip  
	
$(vi)$ If the rank of $(c_{i} u)\Psi$ is $0$, then the result follows by $(ii)$. So let us suppose that $(c_{i} u)\Psi$ has rank $1$. 
Then $\im(c_{i}\Psi)\cap \dom(u\Psi)=\{2i-r\}$ for some $i\leq r \leq \min\{2i-1 ,n\}$. 
Hence, 
$$u\Phi = \left\{ \begin{array}{ll}
	(a_{1}\cdots a_{2i-r-1} a_{2i-r+1} \cdots a_{i} u)\Phi &\mbox{if } i\leq \left\lfloor \frac{n+1}{2} \right\rfloor \vspace{1mm}\\
	(a_{2i-n}\cdots a_{2i-r-1} a_{2i-r+1} \cdots a_i u)\Phi  &\mbox{if } i> \left\lfloor \frac{n+1}{2} \right\rfloor, 
\end{array}\right.$$
and so, by Theorem \ref{t18},  
$$u \,\approx\,  \left\{ \begin{array}{ll} 
	a_{1}\cdots a_{2i-r-1} a_{2i-r+1} \cdots a_{i} u &\mbox{if } i\leq \left\lfloor \frac{n+1}{2} \right\rfloor \vspace{1mm}\\
	a_{2i-n}\cdots a_{2i-r-1} a_{2i-r+1} \cdots a_{i} u  &\mbox{if } i> \left\lfloor \frac{n+1}{2} \right\rfloor. 
\end{array}\right.$$
Just like above, we have $c_{i}\, \sim_{R_{9}} \, c_{i} a_{i+1}\cdots a_{n}$ and, if $i> \left\lfloor \frac{n+1}{2} \right\rfloor$, we also have $c_{i} \, \sim_{R_{9}} \, c_{i} a_{1} \cdots a_{2i-n-1} a_{i+1}\cdots a_{n}$. Thus, by applying also $(v)$, we have 
$$c_{i} u\, \approx\, c_{i} a_{1}\cdots a_{2i-r-1} a_{2i-r+1} \cdots a_{n} u\, \approx\, vu,$$
for some $v\in A^{*}$, as required. \qed
\end{proof}

\begin{lemma}\label{l3predp} 
Let $3\leq j\leq n-1$ and $k,t\geq 0$ be such that  $j-t\geq 3$ and $k-2t\geq 0$. Then $c_{j} b^{k} \approx c_{j-t} b^{k-2t} v_{t}$ for some $v_{t}\in \{a_{1},\ldots, a_{n}\}^{*}$. 
\end{lemma}

\begin{proof}
First, observe that 
\begin{equation} \label{e4}
	a_{i} b^{k}\, \approx\, \left\{\begin{array}{cl} 
		b^{k}          & \mbox{if } i\leq k \qquad (\mbox{by Lemma \ref{l17} $(i)$})\\
		b^{k}  a_{i-k} & \mbox{if } i>k \qquad (\mbox{by } R_{3}). 
	\end{array}\right. 
\end{equation}
Now since the result is trivial for $t=0$, with $v_{t}=\emptyset$, we will proceed by induction on $t\geq1$. Suppose that $j-1\geq 3$ and $k-2\geq0$. Then, by (\ref{e4}), $a_{j-1} b^{k-2}\, \approx\, b^{k-2} v_{1}$, for some $v_{1}\in \{a_{1}, \ldots, a_{n}\}^{*}$, and so 
\begin{equation}\label{e5}
	c_{j} b^{k}\, \sim_{R_{13}}\, c_{j-1} a_{j-1} b^{k-2}\, \approx\, c_{j-1} b^{k-2}v_{1}. 
\end{equation}
Next suppose that $j-t\geq 3$, $k-2t\geq 0$ and $c_{j} b^{k}\, \approx\, c_{j-(t-1)} b^{k-2(t-1)} v_{t-1}$, for some $v_{t-1}\in \{a_{1}, \ldots, a_{n}\}^{*}$. Similar to (\ref{e5}), by taking $j-t+1$ instead of $j$ and $k-2t+2$ instead of $k$, we have $c_{j-(t-1)} b^{k-2(t-1)} \,\approx \,c_{j-t} b^{k-2t} v'_{1}$,  for some $v'_{1}\in \{a_{1}, \ldots, a_{n}\}^{*}$, whence $c_{j} b^{k} \,\approx\, c_{j-t} b^{k-2t} v_{t}$, with $v_{t}= v'_{1} v_{t-1}\in \{a_{1}, \ldots, a_{n}\}^{*}$, as required. \qed
\end{proof}

\begin{lemma}\label{l4predp} 
Let $2\leq i\leq j\leq n-1$ and $u,v\in A^{*}$ be such that $(c_{i} u)\Psi = (c_{j} v)\Psi$. Then $c_{i}u\, \approx\, c_{j} v$. 
\end{lemma}

\begin{proof}
First, we suppose that $i=j$. Then we have $c_{i}\, \sim_{R_{9}} \, c_{i} a_{i+1}\cdots a_{n}$ and, if $i> \left\lfloor \frac{n+1}{2} \right\rfloor$, we also have $c_{i} \, \sim_{R_{9}} \, c_{i} a_{1} \cdots a_{2i-n-1} a_{i+1}\cdots a_{n}$. Let $u'= a_{i+1}\cdots a_{n}u$ and $v'= a_{i+1}\cdots a_{n}v$, if $2\leq i\leq \left\lfloor \frac{n+1}{2} \right\rfloor$, and $u'=a_{1} \cdots a_{2i-n-1} a_{i+1}\cdots a_{n}u$ and $v'=a_{1} \cdots a_{2i-n-1} a_{i+1}\cdots a_{n}v$, if $\lfloor \frac{n+1}{2} \rfloor <i\leq n-1$. Then $c_{i} u'\, \approx\, c_{i} u$,\, $c_{i} v'\, \approx\, c_{i} v$, $\dom(u'\Psi) \subseteq \im(c_{i}\Psi)$ and $\dom(v'\Psi) \subseteq \im(c_{i}\Psi)$. Since $(c_{i} u)\Psi = (c_{j} v)\Psi$, we also get $(c_{i} u')\Psi = (c_{j} v')\Psi$, whence $\dom((c_{i} u')\Psi) = \dom((c_{j} v')\Psi)$, and so  
$$\dom(u'\Psi) =\dom(u'\Psi) \cap \im(c_{i}\Psi) =\dom(v'\Psi) \cap \im(c_{i}\Psi) =\dom(v'\Psi).$$
On the other hand, if $y\in  \dom(u'\Psi) =\dom(v'\Psi)$, then $y=(x)(c_{i}\Psi)$ for some $x\in \dom(c_{i}\Psi)$, and so 
\begin{eqnarray*}
	(y)(u'\Psi) &=& ((x)(c_i\Psi))(u'\Psi) = (x)(c_iu')\Psi = (x)(c_iv')\Psi \\
	&=& ((x)(c_i\Psi))(v'\Psi) =(y)(v'\Psi). 
\end{eqnarray*}
Hence, $u'\Psi=v'\Psi$ and since $u',v'\in A^*$, we also have $u'\Phi=v'\Phi$, and so, by Theorem \ref{t18}, we can conclude that $u'\, \approx\, v'$, from which follows 
$$c_{i}u\, \approx\, c_{i}u'\, \approx\, c_{i}v'\, \approx\, c_{i}v.$$
	
Next we let $i=2$ and $j=3$. As in the previous case, by $R_{9}$, we have $c_{2}\, \approx\, c_{2} a_{3} \cdots a_{n}$ and $c_{3}\, \approx\, c_{3} a_{4} \cdots a_{n}$. Let $u'=a_{3} \cdots a_{n} u$ and $v'=a_{4} \cdots a_{n} v$. Then $c_{2} u'\, \approx\, c_{2} u$,  $c_{3} v'\, \approx\, c_{3} v$, $\dom(u'\Psi)\subseteq\{1,2\}$, $\dom(v'\Psi) \subseteq\{1,2,3\}$ and $(c_{2}u')\Psi =(c_{2}u)\Psi =(c_{3}v)\Psi =(c_{3}v')\Psi$. Since 
$$\dom((c_{2}u')\Psi =\dom((c_{3}v')\Psi) \subseteq \dom(c_{2}\Psi) \cap \dom(c_{3}\Psi)=\{3\},$$ 
it follows that either $u'\Psi= 0_{n} =v'\Psi$ or 
$u'\Psi= \left( \begin{smallmatrix} 
	1\\1 
\end{smallmatrix} \right)$ and 
$v'\Psi= \left( \begin{smallmatrix} 
	3\\1 
\end{smallmatrix} \right)$.  
If $u'\Psi=0_{n}=v'\Psi$, then $(c_{2}u')\Psi=0_{n}=(c_{3}v')\Psi$, and so, by Lemma \ref{l2predp}, 
$$c_{2} u\, \approx\, c_{2} u'\, \approx\, a_{1} \cdots a_{n} \,\approx\, c_{3} v' \,\approx\, c_{3} v. $$
Now suppose that 
$u'\Psi=\left( \begin{smallmatrix}
	1\\1
\end{smallmatrix} \right)$ and 
$v'\Psi=\left( \begin{smallmatrix}
	3\\1
\end{smallmatrix} \right)$. Then, by Theorem \ref{t18}, we have $u'\, \approx\, a_{2} \cdots a_{n}$ and $v'\, \approx\, a_{4}\cdots a_{n} b^{2}$, whence 
\begin{eqnarray*}
	c_{2} u &\approx& c_{2} u'\, \approx\, c_{2} a_{2}\cdots a_{n}\, \sim_{R_{2}}\, c_{2} a_{3}\cdots a_{n} a_{2}\, \sim_{R_{9}}\, c_{2} a_{2}\\
	&\approx& c_{3} a_{4} \cdots a_{n} b^{2} \qquad (\mbox{by Lemma \ref{l1predp}})\\
	&\approx& c_{3} v'\, \approx\, c_{3} v. 
\end{eqnarray*} 
	
Finally, we suppose that $3\leq i < j$. If $(c_{i} u)\Psi =0_{n}$, then $(c_j v)\Psi =0_{n}$, and so, by Lemma \ref{l2predp}, $c_{i} u\, \approx\, a_{1}\cdots a_{n}\, \approx\, c_{j} v$. So suppose that $(c_{i}u)\Psi \neq 0_{n}$. Let us begin by noticing that, just as in (\ref{e4}), using the relations $R_{3}$ and Lemma  \ref{l17} $(i)$, we can obtain $v\, \approx\, b^{k} v_{0}$ for some $k\geq 0$ and $v_{0}\in \{a_{1}, \ldots, a_{n}\}^{*}$. 
	
By contradiction, suppose that $k<2(j-i)$. Let $t= \left\lfloor \frac{k}{2} \right\rfloor$. Then $j-t>j-(j-i) =i\geq 3$ and $0\leq k-2t\leq 1$, and so, by Lemma  \ref{l3predp}, $c_{j} b^{k}\, \approx\, c_{j-t} b^{k-2t} v_{t}$ for some $v_{t}\in\{a_1,\ldots,a_n\}^{*}$, whence 
$$(c_{i}u)\Psi = (c_{j} v)\Psi = (c_{j} b^{k} v_{0})\Psi =  (c_{j-t} b^{k-2t} v_{t} v_{0})\Psi $$
Let $x\in\dom((c_i u)\Psi)$. Then $(x)(c_i u)\Psi=(2i-x)(u\Psi) \leq 2i-x$, and on the other hand, 
\begin{eqnarray*}
	(x)(c_i u)\Psi&=&(x)(c_{j-t} b^{k-2t} v_t v_0)\Psi \\
	&=& (2(j-t)-x)(b^{k-2t} v_t v_0)\Psi\in \{2(j-t) - x, 2(j-t) - x- 1\}.
\end{eqnarray*} 
Since $j-t\geq i+1$, we have $2(j-t)\geq 2i +2$, whence 
$$(x)(c_{i} u)\Psi \geq  2(j-t) - x- 1 \geq 2i-x+1 > 2i-x \geq (x)(c_{i} u)\Psi,$$
which is a contradiction. Therefore, $2(j-i)\leq k$. 
	
Let $t=j-i$. Then $j-t=i \geq 3$ and $k-2t= k-2(j-i) \geq0$. By Lemma \ref{l3predp}, it follows that $c_{j} b^{k}\, \approx\, c_{i} b^{k-2j+2i} v_{t}$, for some $v_{t}\in \{a_{1}, \ldots, a_{n}\}^{*}$, and so $(c_{j} b^{k})\Psi =(c_{i} b^{k-2j+2i} v_{t})\Psi$, which implies that 
$$(c_{i} u)\Psi =(c_{j} v)\Psi =(c_{j}b^{k} v_{0})\Psi =(c_{i} b^{k-2j+2i} v_{t}v_{0})\Psi. $$ 
Thus, by the first part of this proof, we get $c_i u \, \approx\, c_{i} b^{k-2j+2i} v_tv_0$, and so 
$$c_{i} u\, \approx\, c_{i} b^{k-2j+2i} v_{t}v_{0}\, \approx\, c_{j} b^{k} v_{0}\, \approx\, c_{j} v,$$
as required. \qed
\end{proof}

\begin{theorem}\label{predp}
For $n\neq 3$, the monoid $\mathcal{DP}^{-}_{n}$ is defined by the monoid presentation
$$\mathcal{Q}=\langle \, a_{1}, \ldots, a_{n}, b, c_2,\ldots,c_{n-1} \mid R_{1}\-- R_{13}\, \rangle.$$
\end{theorem}

\begin{proof} 
For $n=1,2$,\, $\mathcal{DP}^{-}_{n}=\mathcal{ODP}^{-}_{n}$ and $\mathcal{Q}=\mathcal{P}$, whence the statement coincides with that of Theorem \ref{t18}. Therefore, let us assume that $n\geq 4$. 
	
Let $w,z\in B^{*}$ be such that $w\Psi=z\Psi$. Given that all relations $R_{1} \-- R_{13}$ hold in $\mathcal{DP}^{-}_{n}$ via $\Psi$, the proof is complete by showing that $w\, \approx\, z$. Considering the relations $R_{8}\--R_{11}$, it is easy to show that $w\, \approx\, c_{i_{1}}\cdots c_{i_{k}} u'$, for some $u'\in A^{*}$ and $2\leq i_{1}, \ldots, i_{k}\leq n-1$ ($k\geq 0$). On the other hand, taking into account the relations $R_{5}\--R_{7}$, it is clear that 
$$c_{i_1}\cdots c_{i_k}\, \approx\, \left\{ \begin{array}{ll}
	c_i u'' & \mbox{if $k$ is odd}\\
	u'' & \mbox{if $k$ is even,}
\end{array}\right. $$
for some $2\leq i\leq n-1$ and $u''\in A^{*}$. Therefore, $w\, \approx\, c_{i} u$ or $w\, \approx\, u$, for some $2\leq i\leq n-1$ and $u\in A^{*}$. Similarly, $z\, \approx\, c_{j} v$ or $w\, \approx\, v$, for some $2\leq j\leq n-1$ and $v\in A^{*}$. 
	
If $w\, \approx\, c_{i} u$ and $z\, \approx\, c_{j} v$, then $(c_{i} u)\Psi =w\Psi =z\Psi= (c_{j} v)\Psi$, and so $c_{i} u\, \approx\, c_{j} v$, by Lemma \ref{l4predp}, whence $w\, \approx\, z$. 
	
Next assume that $w\, \approx\, c_{i}u$ and $z\,\approx\, v$. Then $(c_{i}u)\Psi =w\Psi= z\Psi=v\Psi$, and so $(c_{i}u) \Psi$ is both order-preserving and order-reversing, whence $(c_{i}u)\Psi$ has rank at most $1$. Thus, by Lemma \ref{l2predp}, $c_{i}u\, \approx\, v'$ for some $v'\in A^{*}$, and so $v'\Psi= (c_{i}u)\Psi =v\Psi$, which implies that $v'\Phi =v\Phi$, whence $v'\,\approx\, v$, by Theorem \ref{t18}. Therefore,  $w\,\approx\, z$. 
	
If $w\, \approx\, u$ and $z\, \approx\, c_{j} v$, then $w\, \approx\, z$, by applying the same reasoning as in the previous case. 
	
Finally, suppose that $w\, \approx\, u$ and $z\, \approx\, v$. Then $u\Phi = u\Psi = w\Psi = z\Psi = v\Psi = v\Phi$, whence $u\, \approx\, v$, by Theorem \ref{t18}, and so $w\, \approx\, z$, as required. \qed 
\end{proof} 

It is shown in \cite[Theorem 2.17]{AKU3} that 
$$
\lvert DP_{n}^{-} \rvert= 2^{n+1}+2^{n-m+1}+ (n-2m+5)2^{m}-\frac{n^{2}+7n+16}{2},
$$ 
where $m=\left\lfloor \frac{n+1}{2}\right\rfloor$. Using GAP \cite{GAP}, we can verify that, for $n=3$, the monoid defined by the presentation $\mathcal{Q}$ has $14$ elements, while $\lvert \mathcal{DP}^{-}_{3} \rvert=13$. However, by adding the relation $a_2 c_2=b^2$, which holds in $\mathcal{DP}^{-}_{3}$, to the relations of $\mathcal{Q}$, we obtain 
a presentation that defines a monoid with $13$ elements. Therefore, the monoid presentation 
\begin{eqnarray*}
	\langle a_{1}, a_{2}, a_{3}, b, c_{2}\hspace{-2mm} &\mid &\hspace{-2mm} a_{1}^{2}=a_{1},\,  a_{2}^{2}=a_{2},\,  a_{3}^{2}=a_{3},\, a_{1}a_{2}=a_{2}a_{1},\, a_{1}a_{3}=a_{3}a_{1}, \\
	\hspace{-2mm} && \hspace{-2mm} a_{2}a_{3}=a_{3}a_{2},\, ba_{3}=b,\, a_{2}b=ba_{1},\, a_{3}b=ba_{2},\, a_{1}b=b,\, b^{4}=b^{3},\\
	\hspace{-2mm} && \hspace{-2mm} b^{3}=a_{1}a_{2}a_{3},\, a_{1}c_{2}=c_{2},\, c_{2}a_{1}=a_{3}c_{2},\, a_{2}c_{2}= c_{2}a_{2},\, c_{2}a_{3}=c_{2},\\
	\hspace{-2mm} && \hspace{-2mm} c_{2}^{2}=a_{1}a_{3},\, c_{2}^{3}=c_{2}^{2},\, bc_{2}=a_{2}b,\,  c_{2}b=a_{3}b,\, 
	a_{2}c_{2}=b^{2}  \rangle. 
\end{eqnarray*} 
defines $\mathcal{DP}^{-}_{3}$.

\paragraph {Acknowledgment} \hspace*{-2.5pt}
The third author was supported by national funds through the FCT - Funda\c c\~ao para a Ci\^encia e a Tecnologia, I.P., 
under the scope of the projects UID/00297/2025 (https://doi.org/10.54499/UID/00297/2025) and 
UID/PRR/00297/2025 (https://doi.org/10.54499/UID/PRR/00297/2025) (Center for Mathematics and Applications \-- NOVA Math).  




\begin{thebibliography}{99}
	
\bibitem{AKU1} Al-Kharousi, F., Kehinde, R. and Umar, A.: Combinatorial results for certain semigroups of partial isometries of a finite chain. Australasian J. Combin. 58(3), 365--375 (2014).

\bibitem{AKU2} Al-Kharousi, F., Kehinde, R. and Umar, A.: On the semigroup of partial isometries of a finite chain. Commun. Algebra \textbf{44}(2), 639--647 (2016).
https://doi.org/10.1080/00927872.2014.984838

\bibitem{AKU3} Al-Kharousi, F., Kehinde, R. and Umar, A.: Combinatorial results for certain semigroups of order-decreasing partial isometries of a finite chain. Journal of Algebra, Number Theory: Advances and Applications \textbf{19}(1), 29--49 (2018). 

\bibitem{AAIK} Ay\i k, G., Ay\i k, H., Dimitrova I. and Koppitz, J.: On certain semigroups of finite monotone and order-decreasing partial transformations. Asian-Eur. J. Math (online published).

\bibitem{BYA} Bugay, L., Ya\u{g}c\i, M.  and Ay\i k, H.: The ranks of certain semigroups of partial isometries. Semigroup Forum \textbf{97}(2), 214--222 (2018).  
https://doi.org/10.1007/s00233-018-9954-1

\bibitem{FP} Fernandes, V.H. and Paulista, T.: On the monoid of partial isometries of a cycle graph. Turk J Math \textbf{47}(6), 1746--1760 (2023). 
https://doi.org/10.55730/1300-0098.3460	

\bibitem{FQ1} Fernandes, V.H. and Quinteiro, T.M.: Presentations for monoids of finite partial isometries. Semigroup Forum \textbf{93}(1), 97--110 (2016). 
https://doi.org/10.1007/s00233-015-9759-4
	
\bibitem{FQ2} Fernandes, V.H. and Quinteiro, T.M.: On the monoid of partial isometries of a finite star graph. Commun. Algebra \textbf{51}(3), 1028--1048 (2023). 
https://doi.org/10.1080/00927872.2022.2121404

\bibitem{GM} Ganyushkin, O. and Mazorchuk, V.: Classical Finite Transformation Semigroups. Springer-Verlag, London (2009).


\bibitem{GKP} Graham, R.L., Knuth, D.E. and Patashnik, O.: Concrete Mathematics: A Foundation for Computer Science, 2nd edition, Addison-Wesley (1994). 

\bibitem{H} Howie, J.M.: Fundamentals of Semigroup Theory. Oxford University Press, New York  (1995).


\bibitem{KMU} Kehinde, R., Makanjuola S.O. and Umar, A.: On the semigroup of order-decreasing partial isometries of a finite chain. Journal of Algebra, Number Theory: Advances and Applications \textbf{19}(1), 51--65 (2018).


\bibitem{Lallement:1979} Lallement, G.: Semigroups and Combinatorial Applications, John Wiley \& Sons, New York (1979).

\bibitem{Ruskuc:1995} Ru\v skuc, N.: Semigroup Presentations, Ph.D. Thesis, University of St-Andrews, 1995.

\bibitem{SS} Sommanee, W. and Sanwong, J.: On the semigroup of order-preserving partial isometries of a finite chain with restricted range. Semigroup Forum \textbf{104}, 166--179 (2022). 
https://doi.org/10.1007/s00233-021-10237-3 

\bibitem{GAP} The GAP Group. (2025). GAP \-- Groups, Algorithms, and Programming, Version 4.15.1. https://www.gap-system.org 


\end{thebibliography}
\end{document}